# An Iterative Procedure for Optimal Control of Bilinear Systems

H. Ramezanpour[1], S. Setayeshi[1], H. Arabalibeik[2], A. Jajrami[3]

[1]Amirkabir University of Technology, Tehran, Iran

`h.ramezanpour@aut.ac.ir, setayesh@aut.ac.ir`

[2]Research Center for Science and Technology in Medicine (RCSTIM), Tehran University of Medical Sciences, Tehran, Iran.

`arabalibeik@tums.ac.ir`

[3]Ferdowsi University of Mashhad, Mashhad, Iran

`jajarmi@stu-mail.um.ac.ir`

*ABSTRACT*

*This paper presents a new and straightforward procedure for solving bilinear quadratic optimal control problem. In this method, first the original optimal control problem is transformed into a nonlinear two-point boundary value problem (TPBVP) via the Pontryagin's maximum principle. Then, the nonlinear TPBVP is transformed into a sequence of linear time-invariant TPBVPs using the homotopy perturbation method (HPM) and introducing a convex homotopy in topologic space. Solving the latter problems through an iterative process yields the optimal control law and optimal trajectory in the form of infinite series. Finally, sufficient condition for convergence of these series is proved by a theorem. Simplicity and efficiency of the proposed method is shown through an illustrative example.*

## 1.INTRODUCTION

Bilinear systems are a special class of nonlinear systems, in which nonlinear terms are constructed by multiplication of control vector and state vector. An overview of the available control strategies for bilinear systems can be found in [1]-[5]. Besides, optimal control is one of the most active subjects in the control theory. Theory and application of optimal control have been widely used in different fields such as aircraft systems [6], robotic [7], biomedicine [8], etc. However, optimal control of nonlinear systems is a challenging task which has been studied extensively for decades.

For solving nonlinear optimal control problems, indirect methods often lead to solving a Hamilton–Jacobi–Bellman (HJB) equation [9] or a nonlinear two-point boundary value (TPBVP) [10]. In general, an analytic solution does not exist for the HJB equation or nonlinear TPBVP, with the exception of the simplest cases. This has inspired researchers to present some approaches to obtain an approximate solution for the above-mentioned problems as well as obtaining an approximate optimal control for nonlinear systems. An excellent literature review on the methods for approximating the solution of HJB equation is provided in [11]. Besides, to approximate the solution of nonlinear TPBVPs, successive approximation approach (SAA) [12] and sensitivity approach [13] have been introduced recently. In those, a sequence of nonhomogeneous linear time-varying TPBVP's is solved instead of directly solving the nonlinear TPBVP derived from the maximum principle. However, solving time-varying equations is much more difficult than solving time-invariant ones.

The homotopy perturbation method (HPM) was initially proposed by the Chinese mathematician J. H. He [14-15]. This method has been widely used to solve nonlinear problems in different fields [16-18]. In contrast to the perturbation method [19], the HPM is independent upon small/large physical parameters in system model. However, like the other traditional non-perturbation methods, such as the Lyapunov's artificial small parameter method [20] and Adomian's decomposition method [21], uniform convergence of the solution series obtained via the HPM cannot be ensured.

The main objective of this paper is to develop an optimal control design algorithm for a class of bilinear systems with a separated linear part. The algorithm is based on the HPM, where an iterative process is proposed to find a solution sequence of the nonlinear TPBVP, derived from the maximum principle. In this process, the nonlinear terms are considered as known additional disturbances so that the problem is transformed into a sequence of nonhomogeneous linear time-invariant TPBVPs. The approach presented reduces the computational time and avoids the trouble of directly solving the nonlinear TPBVP or the HJB equation. A simulation example is employed to verify the validity of the suggested algorithm.

## 2. PROBLEM STATEMENT

Consider a bilinear system of the form:

$$\begin{cases} \dot{x}(t) = Ax(t) + Bu(t) + \left( \sum_{j=1}^{n} x_j N_j \right) u(t) \\ x(t_0) = x_0 \end{cases} \quad (1)$$

where $x \in R^n$ and $u \in R^m$ are the state and control vectors respectively, $A$, $B$ and $N_j$ are constant matrices of appropriate dimensions with $N_j \in R^{n \times m}$, and $x_0 \in R^n$ is the initial state vector. The objective is to find the optimal control law $u^*(t)$ that minimizes the following quadratic performance index (QPI):

$$J = \frac{1}{2} x^T(t_f) Q_f x(t_f) + \frac{1}{2} \int_{t_0}^{t_f} \left( x^T(t) Q x(t) + u^T(t) R u(t) \right) dt \quad (2)$$

subject to the system (1) where $Q$ and $Q_f$ are symmetric positive semi-definite $n \times n$ matrices, and $R$ is a symmetric positive definite $m \times m$ matrix.

According to the Pontryagin's maximum principle, the optimality condition is obtained as the following nonlinear TPBVP:

$$\begin{cases} \dot{x}(t) = Ax(t) - BR^{-1}B^T \lambda(t) + \phi(x(t), \lambda(t)) \\ \dot{\lambda}(t) = -Qx(t) - A^T \lambda(t) + \psi(x(t), \lambda(t)) \\ x(t_0) = x_0 \\ \lambda(t_f) = Q_f x(t_f) \end{cases} \quad (3)$$

in which:

$$\begin{cases} \phi(x(t), \lambda(t)) = -BR^{-1}\left(\sum_{j=1}^{n} x_j N_j\right)^T \lambda(t) + \left(\sum_{j=1}^{n} x_j N_j\right) w(x(t), \lambda(t)) \\ \psi(x(t), \lambda(t)) = -\frac{\partial}{\partial x}\left(\frac{1}{2} w^T(x(t), \lambda(t)) R w(x(t), \lambda(t)) + \lambda^T \left(B + \sum_{j=1}^{n} x_j N_j\right) z(x)\right) \\ z(x) = -R^{-1}\left(B + \sum_{j=1}^{n} x_j N_j\right)\left(B + \sum_{j=1}^{n} x_j N_j\right)^T \\ w(x(t), \lambda(t)) = z(x)\lambda(t) \end{cases} \quad (4)$$

and $\lambda \in R^n$ is the co-state vector. Also the optimal control law is given by:

$$u^*(t) = -w(x(t), \lambda(t)) \quad (5)$$

Unfortunately, it is difficult to solve problem (3) even numerically. To overcome this difficulty, we will use the HPM in the next section, which transforms problem (3) into a sequence of linear time-invariant TPBVPs. Solving the latter problems recursively and intercepting frontal $N$ terms of the series solution, a suboptimal solution is obtained.

## 3. PROPOSED METHOD

Let us define the operators $F_1[x(t), \lambda(t)]$ and $F_2[x(t), \lambda(t)]$ as follows:

$$F_1[x(t), \lambda(t)] \stackrel{\Delta}{=} \dot{x}(t) - Ax(t) + BR^{-1}B^T \lambda(t) - \phi(x(t), \lambda(t)) \quad (6)$$

$$F_2[x(t), \lambda(t)] \stackrel{\Delta}{=} \dot{\lambda}(t) + Qx(t) + A^T \lambda(t) - \psi(x(t), \lambda(t)). \quad (7)$$

From (3) it is obvious that:

$$\begin{cases} F_1[x(t), \lambda(t)] = 0 \\ F_2[x(t), \lambda(t)] = 0 \end{cases} \quad (8)$$

The operators $F_1$ and $F_2$ can generally be divided into two parts, a linear part and a nonlinear part. Therefore, we can write:

$$\begin{cases} F_1[x(t),\lambda(t)] = L_1[x(t),\lambda(t)] + N_1[x(t),\lambda(t)] \\ F_2[x(t),\lambda(t)] = L_2[x(t),\lambda(t)] + N_2[x(t),\lambda(t)] \end{cases} \quad (9)$$

where $L_i$ and $N_i$ are the linear and nonlinear parts of $F_i$ for $i=1,2$ respectively. Now, we construct two homotopies in topologic space for (9), $\tilde{x}(t,p):[t_0,t_f]\times[0,1]\longrightarrow R$ and $\tilde{\lambda}(t,p):[t_0,t_f]\times[0,1]\longrightarrow R$, which satisfy the following equation:

$$\begin{cases} L_1[(\tilde{x}(t,p),\tilde{\lambda}(t,p))-(x_{ini}(t),\lambda_{ini}(t))] + pL_1[x_{ini}(t),\lambda_{ini}(t)] + pN_1[\tilde{x}(t,p),\tilde{\lambda}(t,p)] = 0 \\ L_2[(\tilde{x}(t,p),\tilde{\lambda}(t,p))-(x_{ini}(t),\lambda_{ini}(t))] + pL_2[x_{ini}(t),\lambda_{ini}(t)] + pN_2[\tilde{x}(t,p),\tilde{\lambda}(t,p)] = 0 \end{cases} \quad (10)$$

with the boundary conditions:

$$\tilde{x}(t_0,p) = x_0, \quad \tilde{\lambda}(t_f,p) = Qx(t_f) \quad (11)$$

where $p \in [0,1]$ is an embedding parameter which is called homotopy parameter, $\lambda_{ini}(t)$ and $x_{ini}(t)$ are the initial guesses for the solution of (3), i.e. $\lambda(t)$ and $x(t)$ respectively. Setting $p=0$ and $p=1$ in (10) yields:

$$p=0 \Rightarrow \begin{cases} L_1[(\tilde{x}(t,0),\tilde{\lambda}(t,0))-(x_{ini}(t),\lambda_{ini}(t))]=0 \\ L_2[(\tilde{x}(t,0),\tilde{\lambda}(t,0))-(x_{ini}(t),\lambda_{ini}(t))]=0 \end{cases} \Rightarrow \begin{cases} \tilde{x}(t,0) = x_{ini}(t) \\ \tilde{\lambda}(t,0) = \lambda_{ini}(t) \end{cases} \quad (12)$$

$$p=1 \Rightarrow \begin{cases} F_1[\tilde{x}(t,1),\tilde{\lambda}(t,1)] = 0 \\ F_2[\tilde{x}(t,1),\tilde{\lambda}(t,1)] = 0 \end{cases} \quad (13)$$

Therefore, if the homotopy parameter $p$ changes from zero to unity, $\tilde{x}(t,p)$ and $\tilde{\lambda}(t,p)$ change from the initial guesses $x_{ini}(t)$ and $\lambda_{ini}(t)$ to the exact solution of (13). In topology we call it deformation. Obviously, when $p=1$, TPBVP (10)-(11) is equivalent to the nonlinear TPBVP (3).

### 3.1 Discussion on $L$ and Initial guesses

It should be mentioned that we have a great freedom to choose linear operators and initial guesses. Once one chooses these parts, the homotopy equation is completely determined, because the remaining part is actually the original equation (see (10)) and we have less freedom to change it. Here we discuss some general rules that should be noted in choosing linear operators and initial guesses.

According to the homotopy perturbation procedure, $L$ should be easy to handle and closely related to the original equation. We mean that it must be chosen in such a way that one has no difficulty in subsequently solving systems of resulting equations. Strictly speaking, in constructing $L$, it is better to use some part of the original equation. A suitable choice of linear operators can be as follows:

$$\begin{cases} L_1[\tilde{x}(t,p), \tilde{\lambda}(t,p)] \overset{\Delta}{=} \dfrac{\partial \tilde{x}(t,p)}{\partial t} - A\tilde{x}(t,p) + BR^{-1}B^T\tilde{\lambda}(t,p) \\ L_2[\tilde{x}(t,p), \tilde{\lambda}(t,p)] \overset{\Delta}{=} \dfrac{\partial \tilde{\lambda}(t,p)}{\partial t} + Q\tilde{x}(t,p) + A^T\tilde{\lambda}(t,p) \end{cases} \quad (14)$$

There is no universal technique for choosing the initial guess in most iterative methods; but from previous works done on HPM and our own experiences, we can conclude that the initial guess should be obtained from the original equation and reduce complexity of the resulting equations. For example, it can be chosen to be the solution to some part of the original equation, or it can be chosen from boundary conditions. One suitable choice can be the solution of the following linear TPBVP, which comes from the original TPBVP (3):

$$\begin{cases} L_1[x_{ini}(t), \lambda_{ini}(t)] = 0 \\ L_2[x_{ini}(t), \lambda_{ini}(t)] = 0 \\ x_{ini}(t_0) = x_0, \lambda_{ini}(t_f) = Qx(t_f). \end{cases} \quad (5)$$

Now the following theorem is presented, which indicate how to use the HPM in practice for handling TPBVP (3):

**Theorem 3.1.** The solution of nonlinear TPBVP (3) can be written as $x(t) = \sum_{n=0}^{\infty} x^{(n)}(t)$ and $\lambda(t) = \sum_{n=0}^{\infty} \lambda^{(n)}(t)$, where the *n*-th order terms $x^{(n)}(t)$ and $\lambda^{(n)}(t)$ for $n \geq 0$ are achieved recursively by solving a sequence of linear time-invariant TPBVPs.

**Proof.** Assume that the embedding parameter $p$ is a small parameter and $\tilde{x}(t,p)$ and $\tilde{\lambda}(t,p)$ are infinitely differentiable with respect to $p$ around $p = 0$. Expanding $\tilde{x}(t,p)$ and $\tilde{\lambda}(t,p)$ as Maclaurin series yields:

$$\begin{cases} \tilde{x}(t,p) = x^{(0)}(t) + x^{(1)}(t)p + x^{(2)}(t)p^2 + \cdots \\ \tilde{\lambda}(t,p) = \lambda^{(0)}(t) + \lambda^{(1)}(t)p + \lambda^{(2)}(t)p^2 + \cdots \end{cases} \quad (16)$$

where $x^{(n)}(t) = \dfrac{1}{n!}\dfrac{\partial^n \tilde{x}(t,p)}{\partial p^n}\bigg|_{p=0}$ and $\lambda^{(n)}(t) = \dfrac{1}{n!}\dfrac{\partial^n \tilde{\lambda}(t,p)}{\partial p^n}\bigg|_{p=0}$. Substituting (15) in (10), rearranging with respect the order of $p$, and equating terms with the same order of $p$ on each side we obtain:

$$p^0 : \begin{cases} L_1[(x^{(0)}(t), \lambda^{(0)}(t))] - L_1[x_{ini}(t), \lambda_{ini}(t)] = 0 \\ L_2[(x^{(0)}(t), \lambda^{(0)}(t))] - L_2[x_{ini}(t), \lambda_{ini}(t)] = 0 \\ x^{(0)}(t_0) = x_0, x^{(0)}(t_f) = x_f \end{cases} \quad (17a)$$

$$p^1 : \begin{cases} L_1[x^{(1)}(t), \lambda^{(1)}(t)] + L_1[x_{ini}(t), \lambda_{ini}(t)] + h_1^{(0)}(t) = 0 \\ L_2[x^{(1)}(t), \lambda^{(1)}(t)] + L_2[x_{ini}(t), \lambda_{ini}(t)] + h_2^{(0)}(t) = 0 \\ x^{(1)}(t_0) = 0, x^{(1)}(t_f) = 0 \end{cases} \quad (17b)$$

$$\vdots$$

$$p^n : \begin{cases} L_1(x^{(n)}(t), \lambda^{(n)}(t)) + h_1^{(n-1)}(t) = 0 \\ L_2(x^{(n)}(t), \lambda^{(n)}(t)) + h_2^{(n-1)}(t) = 0 \\ x^{(n)}(t_0) = 0, x^{(n)}(t_f) = 0 \\ n \geq 2 \end{cases} \quad (17c)$$

where nonhomogeneous terms $h_1^{(n-1)}(t)$ and $h_2^{(n-1)}(t)$ are calculated using the information obtained from previous steps as follows:

$$h_1^{(n-1)}(t) = \frac{1}{(n-1)!}\left(\frac{\partial^{n-1}}{\partial p^{n-1}}\left(N_1\big([x^{(0)} + px^{(1)} + p^2 x^{(2)} + ...], [\lambda^{(0)} + p\lambda^{(1)} + p^2\lambda^{(2)} + . \right.\right. \quad (18)$$

$$h_2^{(n-1)}(t) = \frac{1}{(n-1)!}\left(\frac{\partial^{n-1}}{\partial p^{n-1}}\left(N_2\big([x^{(0)} + px^{(1)} + p^2 x^{(2)} + ...], [\lambda^{(0)} + p\lambda^{(1)} + p^2\lambda^{(2)} + . \right.\right. \quad (19)$$

Therefore, at each step, we obtain a nonhomogeneous linear time-invariant TPBVP which should be solved in a recursive manner.

After obtaining $x^{(n)}(t)$ and $\lambda^{(n)}(t)$ for $n \geq 0$, we should set $p = 1$ in (16) to obtain the exact solution of problem (3). Setting $p = 1$ in (16) yields:

$$\begin{cases} x(t) = \tilde{x}(t,1) = x^{(0)}(t) + x^{(1)}(t) + x^{(2)}(t) + \cdots \\ \lambda(t) = \tilde{\lambda}(t,1) = \lambda^{(0)}(t) + \lambda^{(1)}(t) + \lambda^{(2)}(t) + \cdots \end{cases} \quad (20)$$

and the proof is complete.

**Remark 3.1.** It should be noted that series in (20) converge rapidly for most cases; however, convergence rate depends upon the nonlinear operators. As this method is an iterative method, so the Banach's fixed point theorem can be applied for convergence study of the series (20).

**Theorem 3.2.** (Sufficient condition of convergence) Suppose that $N \stackrel{\Delta}{=} \begin{bmatrix} N_1[\tilde{x}(t,p), \tilde{\lambda}(t,p)] \\ N_2[\tilde{x}(t,p), \tilde{\lambda}(t,p)] \end{bmatrix}$ and $z \stackrel{\Delta}{=} \begin{bmatrix} x(t) \\ \lambda(t) \end{bmatrix}$ and $X$ and $Y$ are Banach spaces and $N : X \to Y$ is a contractive nonlinear mapping that is

$$\forall w, w^* \in X; \|N(w) - N(w^*)\| \leq \gamma \|w - w^*\|, \quad 0 < \gamma < 1$$

Then according to Banach's fixed point theorem $N$ has a unique fixed point $z$ that is $N(z) = z$.

The sequence generated by the HPM may be regarded as

$$W_n = N(W_{n-1}) \quad , \quad W_{n-1} = \sum_{i=0}^{n-1} W_i \quad , \quad n = 1,2,3,...$$

And suppose that $W_0 = w_0 \in B_t(w)$ where $B_t(w) = \{w^* \in X \mid \|w^* - w\| < t \in [t_0, t_f]\}$, then we have

(i) $W_n \in B_t(w)$

(ii) $\lim_{n \to \infty} W_n = w$

**Proof.** (i) By inductive approach, for $n = 1$ we have

$$\|W_1 - w\| = \|N(W_0) - N(w)\| \leq \gamma \|w_0 - w\|$$

Assume that $\|W_{n-1} - w\| \leq \gamma^{n-1} \|w_0 - w\|$, as an induction hypothesis, then

$$\|W_n - w\| = \|N(W_{n-1}) - N(w)\| \leq \gamma \|W_{n-1} - w\| \leq \gamma^n \|w_0 - w\|.$$

Using (i) we have

$$\|W_n - w\| \leq \gamma^n \|w_0 - w\| \leq \gamma^n t \leq t \Rightarrow W_n \in B_t(w)$$

(ii) Because of $\|W_n - w\| \leq \gamma^n \|w_0 - w\|$ and $\lim_{n \to \infty} \gamma^n = 0$, $\lim_{n \to \infty} \|W_n - w\| = 0$, that is, $\lim_{n \to \infty} W_n = w$ and proof is complete.

**Remark 3.2.** Substituting (20) in (4), the optimal control law is obtained as follows:

$$u^*(t) = -R^{-1} B^T \sum_{n=0}^{\infty} \lambda^{(n)}(t). \tag{21}$$

## 4. CASE STUDY

The bilinear model of a chemical reactor [5] is given by:

$$A = \begin{bmatrix} \dfrac{13}{6} & \dfrac{5}{12} \\ -\dfrac{50}{3} & -\dfrac{8}{3} \end{bmatrix} \quad , \quad B = \begin{bmatrix} -\dfrac{1}{8} \\ 0 \end{bmatrix} \quad , \quad N_1 = \begin{bmatrix} -1 \\ 0 \end{bmatrix} \quad , \quad N_2 = \begin{bmatrix} 0 \\ 0 \end{bmatrix}$$

The normalized state variable $x_1$ and $x_2$ represent temperature and concentration of the initial product of the chemical reaction, respectively. The normalized scalar control $u$ represents the cooling flow rate in a jacket around the reactor. The objective is to transfer the system in a finite time very closely to the origin. Weighting matrices are chosen as follows:

$$F = \begin{bmatrix} 1000 & 0 \\ 0 & 1000 \end{bmatrix}, \quad Q = \begin{bmatrix} 10 & 0 \\ 0 & 10 \end{bmatrix}, \quad R = 1$$

Initial conditions are $x(0) = x_0 = [0.15 \ 0]^T$ and $t_f = 3$. Simulation results for three iterations are presented in Figs. 1, 2 and 3, which show achieving desired goals.

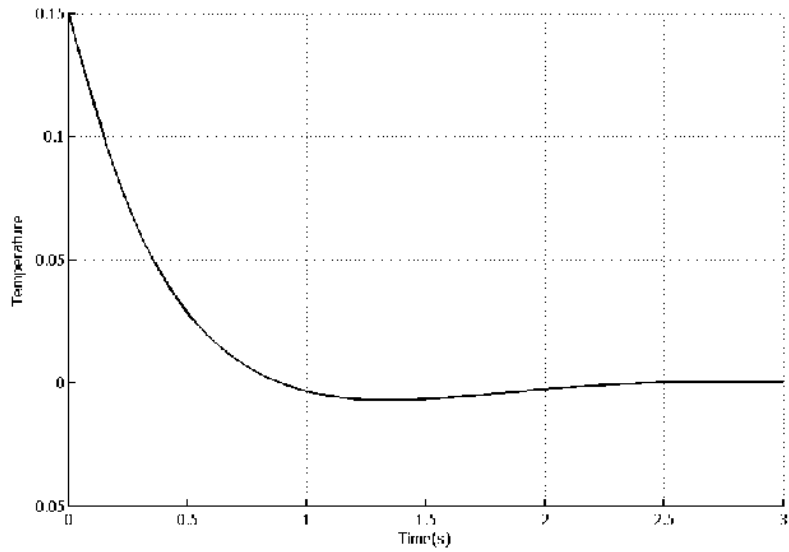

Fig 1. Profiles of temperature

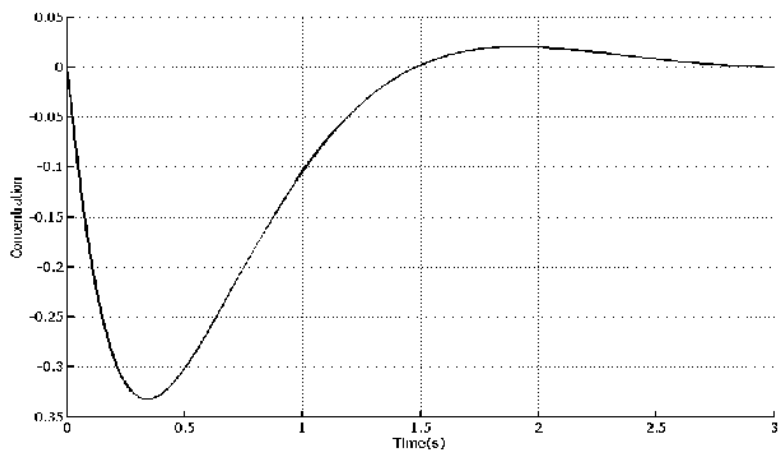

Fig 2. Profiles of Concentration

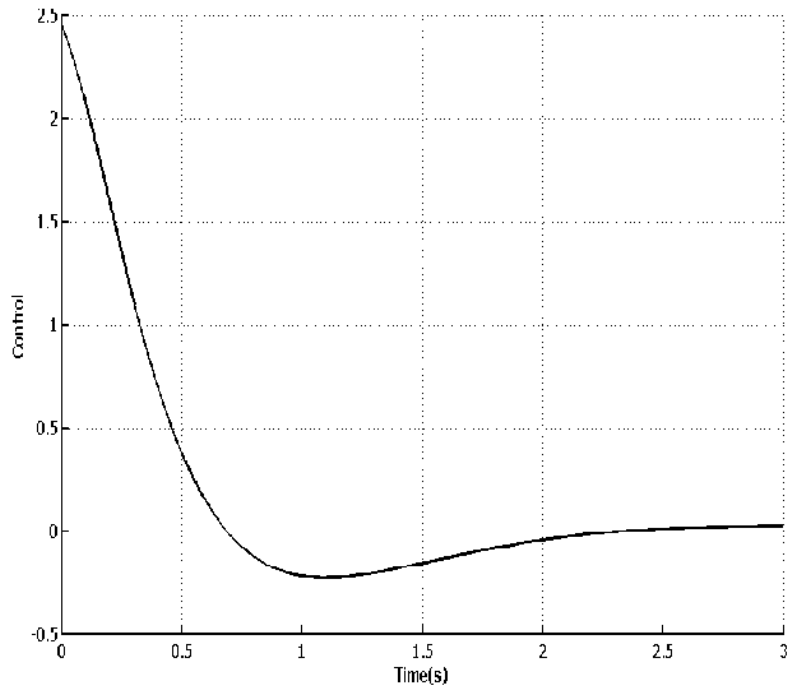

Fig 3. Control Input

## 5. CONCLUSION

In this paper, based on the HPM, an efficient iterative method has been introduced for optimization of bilinear control systems. In this method, by introducing a recursive process, the optimal control law is determined in the form of infinite series with easy-computable terms. The proposed method avoids directly solving the nonlinear TPBVP or the HJB equation. In addition, despite of the successive approximation approach [12] and sensitivity approach [13], it avoids solving a sequence of linear time-varying TPBVPs. It only requires solving a sequence of linear time-invariant TPBVPs. Therefore, in view of computational complexity, the proposed method is more practical than the above-mentioned approximate methods.

## REFRENCES